\documentclass[runningheads]{llncs}
\usepackage[T1]{fontenc}
\usepackage{graphicx}
\usepackage{hyperref}
\usepackage{color}

\usepackage{amsmath}
\usepackage{amssymb}
\usepackage{cleveref}
\usepackage[frozencache]{minted}
\usepackage{tikz-uml}

\newcommand{\CC}{\mathbb C}
\newcommand{\FF}{\mathbb F}
\newcommand{\QQ}{\mathbb Q}
\newcommand{\RR}{\mathbb R}
\newcommand{\ZZ}{\mathbb Z}

\begin{document}
\title{Implementing $p$-adic numbers in Macaulay2 using its foreign function interface and FLINT}

\titlerunning{Implementing $p$-adic numbers in Macaulay2 using its FFI and FLINT}
\author{Douglas A. Torrance\inst{1,2}\orcidID{0000-0003-3999-4973}}
\authorrunning{D. A. Torrance}
%
\institute{Piedmont University, Demorest, GA 30535, USA \and
Georgia Institute of Technology, Atlanta, GA 30332, USA
  \email{dtorrance9@gatech.edu}}

\maketitle              
\begin{abstract}
  Macaulay2 is a computer algebra platform widely used by researchers in algebraic geometry and commutative algebra. Using the ForeignFunctions package, it is possible to make calls from Macaulay2 to dynamic libraries such as FLINT. We demonstrate this by introducing a new Macaulay2 package implementing $p$-adic numbers using FLINT via this interface. We discuss implementation details such as memory allocation, interaction with Macaulay2’s garbage collector, and object-oriented design decisions that mirror the existing implementations of the real and complex number fields in Macaulay2.
\end{abstract}
\section{Introduction}

Macaulay2 \cite{M2} is an important research tool for computational algebraic geometers and commutative algebraists.  It supports computations in a wide variety of rings, including $\ZZ$, $\QQ$, $\RR$, $\CC$, $\FF_{p^k}$, polynomial rings, quotient rings, and fraction fields.  However, historically it has not supported computations in the field $\QQ_p$ of $p$-adic numbers.

Enter FLINT \cite{flint}, a popular C library for number theory, which does support $p$-adic numbers.  Using the Macaulay2 foreign function interface \cite{foreign-functions}, it is possible to make calls to $p$-adic number functions in FLINT from top-level Macaulay2.  We introduce the \textit{Padic} package for Macaulay2, available at \url{https://github.com/d-torrance/macaulay2-padic}, which does this to provide an interface for Macaulay2 users to perform computations in $\QQ_p$.  This package has been distributed with Macaulay2 since version 1.26.05 and provides complete documentation of its features.

This exists as a proof of concept, as potentially $p$-adic support (likely still using FLINT) might one day be added to the Macaulay2 engine and interpreter so that computations in not just $\QQ_p$, but in polynomial rings over $\QQ_p$ and using matrices with entries in $\QQ_p$, are possible.

This paper is organized as follows.  In \Cref{p-adic review}, we provide a brief overview of the $p$-adic numbers, and in \Cref{flint p-adic module} we describe their support in FLINT.  In \Cref{ffi overview}, we describe how Macaulay2 interacts with FLINT at a low level, and in \Cref{object-oriented design}, we look at the interface from a higher level, describing in particular various design decisions.  Finally, we close with a Macaulay2 computation demonstrating Hensel's lemma in \Cref{hensel example}.

\section{A review of $p$-adic numbers}\label{p-adic review}

We open with a crash course on $p$-adic numbers.  The reader is encouraged to reference a textbook such as \cite{gouvea} for more information.

Suppose $p$ is a prime number.  The field of \textit{$p$-adic numbers} $\QQ_p$ is a set whose elements may be expressed uniquely as formal Laurent series
\begin{equation}\label{p-adic expansion}
  \sum_{n=k}^\infty a_n p^n = a_k p^k + a_{k + 1}p^{k + 1} + \cdots,
\end{equation}
where $k\in\ZZ$ and $a_n\in\{0,\ldots,p-1\}$, together with the usual operations of addition and multiplication in base $p$, but noting that carrying may continue indefinitely.  An enlightening exercise is to verify that $-1 = \sum\limits_{n=0}^\infty(p - 1)\cdot p^n$.

For a given $x\in\QQ_p$, its \textit{$p$-adic valuation}, denoted $\nu_p(x)$, is the smallest exponent that appears with a nonzero coefficient in the expansion given in \Cref{p-adic expansion}, e.g., $\nu_3(1\cdot 3^{-2} + 2\cdot 3^{-1})=-2$.  In particular, $\nu_p(0) = \infty$.

These series are not convergent in general using the usual absolute value, but they do converge using the $p$-adic absolute value
\begin{equation*}
|x|_p = p^{-\nu_p(x)}.
\end{equation*}
In fact, $\QQ_p$ is the completion of $\QQ$ under $|\cdot|_p$, just as $\RR$ is the completion of $\QQ$ under the usual absolute value.

The set of all $x\in\QQ_p$ with $\nu_p(x)\geq 0$ forms a ring, known as the \textit{$p$-adic integers} and denoted $\ZZ_p$.  Furthermore, the multiplicative group $\ZZ_p^\times$ of the units of $\ZZ_p$ is precisely the $u\in\QQ_p$ for which $\nu_p(u) = 0$.  It follows that for any nonzero $x\in\QQ_p$ with $\nu_p(x)=\nu$,
\begin{equation}\label{in terms of unit}
  x = \sum_{n=\nu}^\infty a_n p^n = \left(\sum_{n=0}^\infty a_{n + \nu}p^n\right)p^\nu = up^\nu,
\end{equation}
where $u=\sum\limits_{n=0}^\infty a_{n + \nu}p^n\in\ZZ_p^\times$.

The $p$-adic exponential and logarithmic functions are defined using their usual power series expansions provided that they converge, i.e., for $x\in\QQ_p$ with $|x|_p < p^{-\frac{1}{p-1}}$,
\begin{equation*}
  \exp_p x = \sum_{n=0}^\infty\frac{x^n}{n!}
\end{equation*}
and for $x\in\QQ_p$ with $|x-1|_p < 1$,
\begin{equation*}
  \log_p x = \sum_{n=1}^\infty\frac{(-1)^{n-1}(x-1)^n}{n}.
\end{equation*}

Finally, for every $x\in\ZZ_p^\times$ there exists a unique $t\in\ZZ_p^\times$, known as the \textit{Teichm\"uller lift} of $x$, for which $x\equiv t\pmod p$ and $t^p = t$.  Note that $t$ depends only on the congruence class of $x$ modulo $p$.

\section{$p$-adic numbers in FLINT}\label{flint p-adic module}

FLINT (Fast Library for Number Theory) \cite{flint} is a C library implementing a wide variety of number theoretical features.  A number of these features, including ball arithmetic and certain computations over $\ZZ$, $\QQ$, and $\FF_{p^k}$, are already used by the Macaulay2 engine and interpreter.

FLINT's $p$-adic number module was written by Sebastian Pancratz.  To represent $x\in\QQ_p$, the expansion from \Cref{in terms of unit} is used.  However, since a given $u\in\ZZ_p^\times$ may require infinitely many digits in base $p$ to represent, it is truncated up to a given precision $N$.  In other words, the elements $x\in\QQ_p$ that may be represented in FLINT are those of the form $x = up^\nu$ where
\begin{equation}\label{p-adic with precision}
  u = \sum_{n = 0}^{N-1}a_{n+\nu}p^n.
\end{equation}

Since $u$ has finitely many nonzero digits, we certainly have $u\in\ZZ$, and so it is sufficient to represent $x$ using the integers $u$, $\nu$, $N$, and $p$.

FLINT uses two main types in its $p$-adic number module.  A \texttt{padic\_t} object contains $u$, $\nu$, and $N$, and a \texttt{padic\_ctx\_t} object contains information about $p$.  This information is separated so that the same \texttt{padic\_ctx\_t} object may be used for computations involving multiple \texttt{padic\_t} objects.

A variety of functions exist for performing operations on $p$-adic numbers in FLINT.  See \Cref{flint p-adic functions} for a selection.

\begin{table}[htbp!]
  \centering
  \begin{tabular}{|c|c|}
    \hline
    operation & FLINT function \\ \hline\hline
    $-x$ & \texttt{padic\_neg} \\
    $x^{-1}$ & \texttt{padic\_inv} \\
    $x+y$ & \texttt{padic\_add} \\
    $x-y$ & \texttt{padic\_sub} \\
    $x\times y$ & \texttt{padic\_mul} \\
    $x\div y$ & \texttt{padic\_div} \\
    $x^n$ $(n\in\ZZ)$ & \texttt{padic\_pow\_si} \\
    $\sqrt x$ & \texttt{padic\_sqrt} \\
    $\exp_p x$ & \texttt{padic\_exp} \\
    $\log_p x$ & \texttt{padic\_log} \\
    Teichm\"uller lift & \texttt{padic\_teichmuller} \\
    $x = y$ & \texttt{padic\_equal} \\
    $x = 0$ & \texttt{padic\_is\_zero} \\
    \hline
  \end{tabular}
  \vspace{1em}
  \caption{FLINT functions for operations in $\QQ_p$}
  \label{flint p-adic functions}
\end{table}

Furthermore, functions exist for converting back and forth between \texttt{padic\_t} objects and \texttt{fmpz\_t} and \texttt{fmpq\_t} objects, which are FLINT's implementations of integers and rationals, respectively.  In particular, using the natural inclusion map $\ZZ\hookrightarrow\QQ_p$ (resp. $\QQ\hookrightarrow\QQ_p$), we may promote an element of $\ZZ$ (resp. $\QQ$) to $\QQ_p$ using \texttt{padic\_set\_fmpz} (resp. \texttt{padic\_set\_fmpq}), and if possible, we may lift an element of $\QQ_p$ to $\ZZ$ (resp. $\QQ$) using \texttt{padic\_get\_fmpz} (resp. \texttt{padic\_get\_fmpq}).

\section{Foreign function interface and memory allocation}\label{ffi overview}

The \textit{ForeignFunctions} package in Macaulay2 \cite{foreign-functions} allows users to make calls to functions in shared libraries such as FLINT without writing any C code or re-compiling the Macaulay2 binary.  In particular, we may use it to call the $p$-adic functions outlined in \Cref{flint p-adic module}.

Macaulay2 uses the Boehm-Demers-Weiser conservative garbage collector \cite{bdwgc} to allocate and, as needed, deallocate memory.  However, each of the various FLINT types used in the $p$-adic number module has functions for both initializing and clearing the memory used by instances of that type.  For example, \texttt{padic\_init} (or \texttt{padic\_init2} to specify the precision) is used to initialize a \texttt{padic\_t} object and \texttt{padic\_clear} is used to clear it.

Therefore, after allocating the memory for a FLINT object and calling the corresponding initialization function, we must register a finalizer so that when the garbage is collected and the memory is freed, the appropriate FLINT function is also called. See \Cref{garbage collection example}.  Note that when allocating memory for a \texttt{padic\_t} object, we are allocating memory for a struct that fits three \texttt{long} integers -- $u$, $\nu$, and $N$ from \Cref{p-adic with precision}.\footnote{Actually, $u$ is a \texttt{fmpz\_t} object, which is either a \texttt{long} integer or, if needed, a pointer to an arbitrary precision integer.}

\begin{listing}[htbp!]
  \begin{minted}[frame=lines,framesep=1ex]{macaulay2}
needsPackage "ForeignFunctions"
flint = openSharedLibrary "flint"
padicInit = foreignFunction(flint, "padic_init", void, voidstar)
padicClear = foreignFunction(flint, "padic_clear", void, voidstar)
y = getMemory(3 * size long)
padicInit y
registerFinalizer(y, padicClear)
\end{minted}
\caption{Initializing a \texttt{padic\_t} object}
\label{garbage collection example}
\end{listing}

\section{Object-oriented design}\label{object-oriented design}

Just using the foreign function interface to call FLINT functions is not sufficient for doing computations over $\QQ_p$ in Macaulay2.  It is necessary to provide a top-level interface with which users will interact and that wraps the low-level calls to FLINT via the foreign function interface.

Because of the similarities with the real and complex fields in that we are dealing with numbers with potentially infinitely many digits, let us first review how $\RR$ and $\CC$ are implemented in Macaulay2.  Real and complex numbers are instances of the classes \texttt{RR} and \texttt{CC}, respectively.  These classes are both instances of \texttt{InexactFieldFamily} but subclasses of \texttt{InexactNumber}.\footnote{There also exists objects representing the fields $\RR$ and $\CC$ that depend on the floating-point precision used, e.g., \texttt{RR\_53} and \texttt{CC\_53} for double precision.  We have not yet implemented analogous objects for $\QQ_p$.}

\begin{figure}[htbp!]
  \centering
  \begin{tikzpicture}
    \umlemptyclass[y=4]{Type}
    \umlemptyclass[y=2]{RingFamily}
    \umlemptyclass[x=-2]{InexactFieldFamily}
    \umlemptyclass[x=2]{PadicFieldFamily}
    \umlinherit{RingFamily}{Type}
    \umlinherit{InexactFieldFamily}{RingFamily}
    \umlinherit{PadicFieldFamily}{RingFamily}
  \end{tikzpicture}
  \caption{Class diagram for number types}
  \label{type class diagram}
\end{figure}
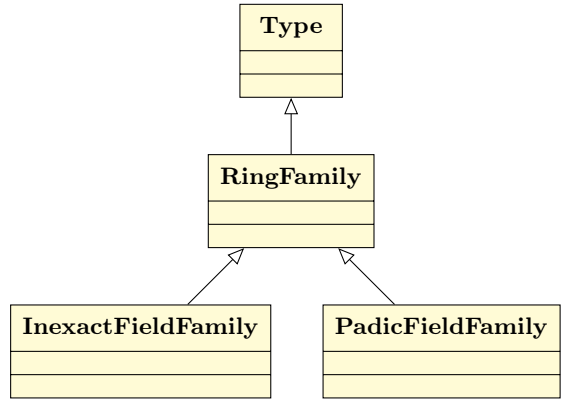

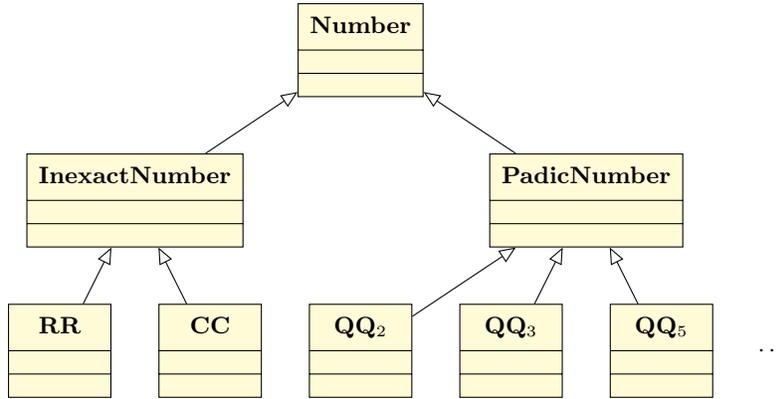
\begin{figure}[htbp!]
  \centering
  \begin{tikzpicture}
    \umlemptyclass[x=1,y=4]{Number}
    \umlemptyclass[x=-2,y=2]{InexactNumber}
    \umlemptyclass[x=4,y=2]{PadicNumber}
    \umlemptyclass[x=-3]{RR}
    \umlemptyclass[x=-1]{CC}
    \umlemptyclass[x=1]{QQ$_2$}
    \umlemptyclass[x=3]{QQ$_3$}
    \umlemptyclass[x=5]{QQ$_5$}
    \umlinherit{InexactNumber}{Number}
    \umlinherit{PadicNumber}{Number}
    \umlinherit{RR}{InexactNumber}
    \umlinherit{CC}{InexactNumber}
    \umlinherit{QQ$_2$}{PadicNumber}
    \umlinherit{QQ$_3$}{PadicNumber}
    \umlinherit{QQ$_5$}{PadicNumber}
    \draw (6.5, 0) node {$\cdots$};
  \end{tikzpicture}
  \caption{Class diagram for number objects}
  \label{object class diagram}
\end{figure}

To parallel this, for a given prime number \texttt{p}, we let the class \texttt{QQ\_p} be an instance of \texttt{PadicFieldFamily} and a subclass of \texttt{PadicNumber}.  See the class diagrams in \Cref{type class diagram} and \Cref{object class diagram}.

Instances of \texttt{QQ\_p} are hash tables containing pointers to the corresponding \texttt{padic\_t} and \texttt{padic\_ctx\_t} objects.  The latter are memoized to avoid duplication.  Then methods are installed on \texttt{PadicNumber} wrapping each of the foreign functions from FLINT, e.g., \texttt{PadicNumber + PadicNumber} calls \texttt{padic\_add}.

The constructor methods for these classes either take a single Macaulay2 number, returning a $p$-adic number with the default precision of 20, or optionally an integer first argument to indicate the precision.  See \Cref{constructor methods}.

\begin{listing}[htbp!]
\begin{minted}[frame=lines,framesep=1ex]{macaulay2}
i1 : needsPackage "Padic";

i2 : QQ_5(-1)

o2 = 4 + 4*5^1 + 4*5^2 + 4*5^3 + 4*5^4 + 4*5^5 + 4*5^6 + 4*5^7 + 4*5^8
     + 4*5^9 + 4*5^10 + 4*5^11 + 4*5^12 + 4*5^13 + 4*5^14 + 4*5^15 +
     4*5^16 + 4*5^17 + 4*5^18 + 4*5^19

o2 : QQ  (of precision 20)
       5

i3 : QQ_5(40, -1)

o3 = 4 + 4*5^1 + 4*5^2 + 4*5^3 + 4*5^4 + 4*5^5 + 4*5^6 + 4*5^7 + 4*5^8
     + 4*5^9 + 4*5^10 + 4*5^11 + 4*5^12 + 4*5^13 + 4*5^14 + 4*5^15 +
     4*5^16 + 4*5^17 + 4*5^18 + 4*5^19 + 4*5^20 + 4*5^21 + 4*5^22 +
     4*5^23 + 4*5^24 + 4*5^25 + 4*5^26 + 4*5^27 + 4*5^28 + 4*5^29 +
     4*5^30 + 4*5^31 + 4*5^32 + 4*5^33 + 4*5^34 + 4*5^35 + 4*5^36 +
     4*5^37 + 4*5^38 + 4*5^39

o3 : QQ  (of precision 40)
       5
\end{minted}
\caption{Constructing $-1\in\QQ_5$ with different precisions}
\label{constructor methods}
\end{listing}
\section{Example: Hensel's lemma}\label{hensel example}

We close with an example, showing how to use our Macaulay2 package to find roots of polynomials over $\ZZ_p$.  We first recall Hensel's lemma, a fundamental result in $p$-adic analysis.  

\begin{lemma}[Hensel \cite{hensel}]
  Suppose $f\in\ZZ_p[x]$.  If there exists $\alpha_1\in\ZZ_p$ for which $f(\alpha_1)\equiv 0\pmod p$ and $f'(\alpha_1)\not\equiv 0\pmod p$, then there exists $\alpha\in\ZZ_p$ satisfying $\alpha_1\equiv\alpha\pmod p$ and $f(\alpha)=0$.

  Furthermore, $\alpha$ may be constructed by defining the recursive sequence given by Newton's method,
  \begin{equation*}
    \alpha_{n+1} = \alpha_n - \frac{f(\alpha_n)}{f'(\alpha_n)}.
  \end{equation*}
Then $f(\alpha_n)\equiv 0\pmod{p^n}$ for all $n$ and $\alpha=\lim\limits_{n\to\infty}\alpha_n$, where the limit is taken with respect to $|\cdot|_p$.
\end{lemma}

As with the Teichm\"uller lift, any choice of $\alpha_1$ in the same congruence class modulo $p$ will work, and in particular, it is natural to choose $\alpha_1\in\{0,\ldots,p-1\}$ so that $\alpha = \alpha_1\cdot p^0 + \cdots$ (see \Cref{p-adic expansion}).

We would like to use Macaulay2 to compute the cube root of 2 in $\ZZ_5$ by finding a root of $x^3-2\in\ZZ_5[x]$.  Since $3^3-2=25\equiv 0\pmod 5$ and $3\cdot 3^2 = 27\not\equiv 0\pmod 5$, we begin with $\alpha_1 = 3$.  We then iterate Newton's method until it converges (up to the $5^{19}$ term, as FLINT's default precision for $p$-adics is $N=20$).  See \Cref{hensel lemma example code}.

\begin{listing}[htbp!]
\begin{minted}[frame=lines,framesep=1ex]{macaulay2}
i1 : needsPackage "Padic";

i2 : newton = x -> x - (x^3 - 2)/(3*x^2);

i3 : alpha = QQ_5 3

o3 = 3

o3 : QQ  (of precision 20)
       5

i4 : while alpha != (alpha = newton alpha) do null

i5 : alpha

o5 = 3 + 2*5^2 + 2*5^3 + 3*5^4 + 1*5^5 + 4*5^6 + 2*5^8 + 3*5^9 +
     4*5^12 + 4*5^14 + 4*5^15 + 3*5^16 + 1*5^17 + 1*5^18 + 2*5^19

o5 : QQ  (of precision 20)
       5

i6 : alpha^3 -- check

o6 = 2

o6 : QQ  (of precision 20)
       5
\end{minted}
\caption{Computing the cube root of 2 in $\ZZ_5$ using Macaulay2}
\label{hensel lemma example code}
\end{listing}

\begin{credits}
  \subsubsection{\ackname} The author thanks the anonymous referees for their helpful comments.  This work was supported by a grant from the Simons Foundation International (SFI-MPS-SSRFA-00012695).

\subsubsection{\discintname}
The author has no competing interests to declare that are
relevant to the content of this article.
\end{credits}
%
%
%
\bibliographystyle{splncs04}
\bibliography{macaulay2-padic}

\end{document}